\pdfoutput=1
\documentclass[11pt]{amsart}
\usepackage[T1]{fontenc}
\usepackage[utf8]{inputenc}
\usepackage{lmodern}
\usepackage[margin=1in]{geometry}
\usepackage{mathtools,amssymb,graphicx,booktabs,microtype,enumitem,array,placeins,xurl}
\usepackage[hidelinks]{hyperref}
\hypersetup{pdftitle={Exact Completion of Locally Few-Distance Sets in Odd Regular Polygons},pdfauthor={Jeonghwan Lee}}
\theoremstyle{plain}
\newtheorem{thm}{Theorem}[section]
\newtheorem{prop}[thm]{Proposition}
\newtheorem{cor}[thm]{Corollary}
\newtheorem{lemma}[thm]{Lemma}
\newtheorem{suptheorem}{Theorem}[section]
\newtheorem{suplemma}[suptheorem]{Lemma}
\newtheorem{supproposition}[suptheorem]{Proposition}
\newtheorem{supcorollary}[suptheorem]{Corollary}
\theoremstyle{remark}
\newtheorem{rem}[thm]{Remark}

\newcommand{\Z}{\mathbb Z}
\newcommand{\sig}{\sigma}

\newcommand{\del}{\mathbin{\triangle}}

\renewcommand{\theequation}{\arabic{equation}}
\setlist{nosep}
\title[Exact completion in odd regular polygons]{Exact Completion of Locally Few-Distance Sets in Odd Regular Polygons}
\author{Jeonghwan Lee}
\address{Korea Advanced Institute of Science and Technology (KAIST), 291 Daehak-ro, Yuseong-gu, Daejeon 34141, Republic of Korea}
\email{jameslee@kaist.ac.kr}
\date{November 27, 2026}
\begin{document}
\begin{abstract}
Let \(A\) be an \(n\)-point subset of an odd regular polygon, let \(k\) be its maximum pinned distance count, and put \(t = 2k - n + 1\). If \(n > \max\{ 6t,(t + 1)(t + 2)/2\}\), then the smallest regular subpolygon containing \(A\) has exactly \(n + t = 2k + 1\) vertices, and \(A\) determines exactly \(k\) distances globally. We prove the underlying completion theorem in every finite abelian group of odd order. The complement in the minimal completion gives an exact obstruction: completion with \(t\) missing points fails precisely when the complement is a sum cover. Representation counts yield explicit bounds for larger local-to-global distance gaps. Finally, we identify the optimal asymptotic completion coefficient with the extremal efficiency of sum covers in odd cyclic groups. Its numerical value remains open.

\end{abstract}
\maketitle
\section{Introduction}\label{introduction}

A set is locally \(k\)-distance if each point determines at most \(k\) distinct nonzero distances to the other points. For \(n\) points on a circle, every distance from a fixed point occurs at most twice, so \(n \leq 2k + 1\). We study the nonnegative defect \(t = 2k - n + 1\) for subsets of odd regular polygons.

We prove that \(n > \max\{ 6t,(t + 1)(t + 2)/2\}\) forces the smallest regular subpolygon containing the set to have exactly \(n + t = 2k + 1\) vertices. Moreover, the global distance count equals the maximum pinned count. A pointwise identity for the missing vertices characterises failure by a sum-cover condition and yields a hierarchy of distance gaps. We then identify the optimal asymptotic completion coefficient with odd cyclic sum-cover efficiency. All proofs of the results stated in this article are included below.

Pinned-distance results for arbitrary convex polygons \cite{dumitrescu,nivasch} do not assume regular ambient positions. Our inverse completion result concerns this more restrictive setting and does not improve their general bounds. Nozaki and Shinohara \cite{nozaki}, building on spherical distance-set theory \cite{delsarte}, study local distance conditions at an equality endpoint. Momihara and Shinohara \cite{momihara,shinohara} obtain regular-polygon containment from sufficiently few global distances on a circle. Here equality of local and global counts is a conclusion, so neither result unconditionally strengthens the other.

The algebraic background includes small sumsets \cite{kneser,tao} and near-maximal three-term progression counts \cite{green-sisask,chase}. Our hypothesis bounds the reflection defect at every centre and determines an exact number of missing points. Sum covers and extremal cyclic bases are established subjects \cite{fitch,haanpaa,jia}; the contribution proposed here is their precise connection to local distance completion. The endpoint of our asymptotic identification is the odd-modulus version of Green's Problem~33 \cite{green}.

Additional finite computations, bounds, and extensions are provided in Sections~S1--S10; none is required for the proofs in Sections~1--8. Throughout, exact completion denotes an exact structural conclusion under sufficient hypotheses; it does not assert an optimal numerical threshold at every defect.

\section{Definitions and the completion theorem}\label{definitions-and-the-completion-theorem}

Let \(\Gamma\) be a finite abelian group of odd order, written additively, and let \(A \subseteq \Gamma\) be nonempty. Write \(n = |A|\). For \(a \in A\), define

\begin{equation}
R_{a}(x) = 2a - x,\quad\quad\sigma_{A}(a) = |A\backslash R_{a}(A)|,\quad\quad t = \max_{a \in A}\sigma_{A}(a).\label{eq:1}
\end{equation}

The map \(R_{a}\) is an involution. Since \(|R_{a}(A)| = |A|\), \(\sigma_{A}(a) = \frac{1}{2}|A \triangle R_{a}(A)|\). Let

\[D = A - A,\quad\quad H = \langle A - A\rangle,\quad\quad C = a_{0} + H\quad(a_{0} \in A).\]

The coset \(C\) is independent of the choice of \(a_{0}\) and is the smallest coset of a subgroup containing \(A\). Indeed, every coset containing \(A\) has underlying subgroup containing \(A - A\) and hence \(H\). Put

\begin{equation}
M = |H|,\quad\quad u = M - n,\quad\quad Q(t) = \max\{ 6t,\frac{(t + 1)(t + 2)}{2}\}.\label{eq:2}
\end{equation}

Thus \(u\) is the number of holes in the minimal affine completion.

\begin{thm}[Exact affine completion]\label{res:2.1}
Let \(A\) be a subset of a finite abelian group of odd order, with the notation above. If \(n > Q(t)\), then

\[A - A = H,\quad\quad|H| = n + t,\quad\quad|C\backslash A| = t.\]

In particular, the completion and its number of holes are determined exactly.

\end{thm}

For the geometric interpretation, identify \(\mathbb{Z}_{m}=\mathbb{Z}/m\mathbb{Z}\) with the vertices \(\exp(2\pi ia/m)\) of the unit regular \(m\)-gon. Define

\[d_{A}(a) = |\{\,|e^{2\pi ia/m} - e^{2\pi ib/m}|:b \in A\backslash\{ a\}\,\}|,\quad\quad k = \max_{a \in A}d_{A}(a).\]

The \emph{global} number of distances is denoted by \(K(A)\).

\begin{cor}[Local-to-global completion]\label{res:2.2}
Let \(m \geq 3\) be odd, let \(A \subseteq \mathbb{Z}_{m}\), and assume \(n = |A| \geq 2\). Set \(t = 2k - n + 1\). If \(n > Q(t)\), then \(A\) is obtained by deleting exactly \(t\) vertices from its unique smallest regular subpolygon, which has \(n + t = 2k + 1\) vertices. Moreover,

\[K(A) = k,\quad\quad 2k + 1 \mid m.\]

\end{cor}

Here and below, a regular subpolygon uses vertices of the given ambient polygon and has the same circumcircle. The word ``smallest'' refers to inclusion, equivalently to the coset \(C\). This does not claim that \(A\) itself is invariant under the rotations or reflections of \(C\).

\section{Reflection defects and subgroup closure}\label{reflection-defects-and-subgroup-closure}

\begin{lemma}[Pinned-distance identity]\label{res:3.1}
For odd \(m\) and every \(a \in A \subseteq \mathbb{Z}_{m}\),

\begin{equation}
2d_{A}(a) = n - 1 + \sigma_{A}(a).\label{eq:3}
\end{equation}

For any finite abelian group of odd order, \(\sigma_{A}(a) \equiv n - 1\ (mod\ 2)\).

\end{lemma}

\begin{proof}
In an odd-order group the unique fixed point of \(R_{a}\) is \(a\). All other orbits have size two. In \(A\backslash\{ a\}\), let \(p\) be the number of complete two-point orbits and \(s\) the number of singleton intersections with these orbits. Then \(n - 1 = 2p + s\) and \(\sigma_{A}(a) = s\), proving the parity assertion. In the cyclic geometric case, two different vertices \(b,c\) have equal distance from \(a\) if and only if \(b + c = 2a\ (mod\ m)\). Thus each of the \(p + s\) orbits that meets \(A\backslash\{ a\}\) contributes exactly one distance. Therefore \(d_{A}(a) = p + s\), which gives (3). \end{proof}

For \(h \in \Gamma\), write \(r_{A}(h) = |A \cap (A + h)|\). We use the metric \(\Delta(X,Y) = \frac{1}{2}|X \triangle Y|\) on finite subsets; it is invariant under bijections and satisfies the triangle inequality.

\begin{lemma}[Uniform defect forces a subgroup difference set]\label{res:3.2}
If \(n > 4t\), then \(D = A - A\) is a subgroup and hence \(D = H\). Furthermore,

\begin{equation}
r_{A}(h) \geq n - 2t\quad(h \in H),\quad\quad|H| \leq 1 + \frac{n(n - 1)}{n - 2t}.\label{eq:4}
\end{equation}

\end{lemma}

\begin{proof}
For \(a,b \in A\), the composition \(R_{a}R_{b}\) is translation by \(2(a - b)\). Consequently,

\[\begin{aligned}
n - r_{A}(2(a - b)) & = \Delta(A,R_{a}R_{b}(A)) \\
 & \leq \Delta(A,R_{a}(A)) + \Delta(R_{a}(A),R_{a}R_{b}(A)) \\
 & = \sigma_{A}(a) + \sigma_{A}(b) \leq 2t.
\end{aligned}\]

Since \(n > 2t\), every element of \(2D\) occurs as a difference of two elements of \(A\), so \(2D \subseteq D\). Doubling is injective in an odd-order group. The finite sets \(2D\) and \(D\) therefore have the same size, and \(2D = D\). It follows that \(r_{A}(h) \geq n - 2t > n/2\) for every \(h \in D\).

Take \(h,k \in D\). The subsets \(A \cap (A + h)\) and \(A \cap (A - k)\) of \(A\) both have more than \(n/2\) elements and hence intersect. If \(z\) lies in their intersection, then \(z - h,z + k \in A\), so \(h + k \in D\). Since \(0 \in D\) and \(D = - D\), this proves that \(D\) is a subgroup. Finally, all \(n^{2}\) ordered differences lie in \(H = D\), and

\[n^{2} = \sum_{h \in H}^{}r_{A}(h) \geq n + (|H| - 1)(n - 2t).\]

This proves (4). \end{proof}

\section{Counting holes and proving exact completion}\label{counting-holes-and-proving-exact-completion}

Translate \(C\) to \(H\) and put \(B = H\backslash A\). Translation preserves the reflection defects. For a subset \(S \subseteq H\), let

\[T_{3}(S) = |\{(a,b) \in S^{2}:2a - b \in S\}|.\]

This counts ordered three-term progressions, including the \(|S|\) constant progressions. Also write \(r_{B + B}(s) = |\{(b,c) \in B^{2}:b + c = s\}|\).

\begin{lemma}[Complement identities]\label{res:4.1}
For every \(A \subseteq H\) as above and every \(a \in A\),

\begin{align}
\sigma_{A}(a) & = u - r_{B + B}(2a),\label{eq:5} \\
\sum_{a \in A}^{}\sigma_{A}(a) & = u(n - u) + T_{3}(B).\label{eq:6}
\end{align}

In particular,

\begin{equation}
t \leq u,\quad\quad u \equiv t\ (mod\ 2),\quad\quad nt \geq u(n - u + 1).\label{eq:7}
\end{equation}

\end{lemma}

\begin{proof}
Since \(R_{a}\) preserves \(H\),

\[A\backslash R_{a}(A) = (H\backslash B) \cap R_{a}(B).\]

The right-hand side has size \(u - |B \cap R_{a}(B)|\), and its subtracted term is \(r_{B + B}(2a)\). This proves (5). The involution interchanging \(b,c\) has no fixed pair in the representations of \(2a\) by \(B + B\): a diagonal pair would require \(b = c = a\), since doubling is injective, whereas \(a \notin B\). Thus \(r_{B + B}(2a)\) is even. It follows that \(\sigma_{A}(a) \leq u\) and \(\sigma_{A}(a) \equiv u\ (mod\ 2)\), giving the first two claims in (7).

Expand \(T_{3}(A)\) using \(1_{A} = 1_{H} - 1_{B}\) in

\[\sum_{a,b \in H}^{}1_{A}(a)1_{A}(b)1_{A}(2a - b).\]

Each single \(1_{B}\) term sums to \(Mu\), each product of two such terms sums to \(u^{2}\), and the triple product sums to \(T_{3}(B)\). For the pair terms, the relevant changes of coordinates are bijective because \(2\) is invertible on \(H\). Hence

\[T_{3}(A) = M^{2} - 3Mu + 3u^{2} - T_{3}(B).\]

On the other hand, \(\sum_{a \in A}^{}\sigma_{A}(a) = n^{2} - T_{3}(A)\). Substitution of \(M = n + u\) gives (6). Since \(T_{3}(B) \geq u\) and every defect is at most \(t\), the last inequality in (7) follows. \end{proof}

\begin{prop}[Quantitative completion before exactness]\label{res:4.2}
If \(n > 6t\), then \(D = H\), \(u < n/2\), and

\begin{equation}
t \leq u \leq \frac{n + 1 - \sqrt{{(n + 1)}^{2} - 4nt}}{2}.\label{eq:8}
\end{equation}

The upper bound can be rounded down to the largest integer congruent to \(t\) modulo \(2\).

\end{prop}

\begin{proof}
Lemma~\ref{res:3.2} gives

\[u = M - n \leq \frac{2nt}{n - 2t} < \frac{n}{2}.\]

The last strict inequality follows from \(n > 6t\), and also holds when \(t = 0\). The function \(f(v) = v(n - v + 1)\) is strictly increasing for \(0 \leq v \leq n/2\). Lemma~\ref{res:4.1} gives \(f(u) \leq nt\). Solving this quadratic inequality on that interval proves (8). Its discriminant is positive because \(n > 6t\). The parity refinement follows from (7). \end{proof}

\begin{proof}[Proof of Theorem~\ref{res:2.1}]
The hypothesis implies \(n > 6t\). Thus \(D = H\) and \(u < n/2\). By Lemma~\ref{res:4.1}, \(u \geq t\) and \(u \equiv t\ (mod\ 2)\). Suppose \(u > t\). Then \(u \geq t + 2\), and monotonicity of \(f(v) = v(n - v + 1)\) up to \(n/2\) yields

\[nt \geq f(u) \geq (t + 2)(n - t - 1).\]

Rearranging gives \(2n \leq (t + 1)(t + 2)\), contradicting \(n > Q(t)\). Therefore \(u = t\), as required. \end{proof}

\begin{proof}[Proof of Corollary~\ref{res:2.2}]
Lemma~\ref{res:3.1} identifies the algebraic defect with \(2k - n + 1\). Theorem~\ref{res:2.1} then gives \(|H| = n + t = 2k + 1\). A subgroup of \(\mathbb{Z}_{m}\) of order \(M\) is the set of multiples of \(m/M\); its cosets are precisely regular \(M\)-gons on the ambient circumcircle. Since \(A - A = H\), all nonzero differences in \(H\) occur. The pairs \(\{ h, - h\}\), \(h \neq 0\), correspond bijectively to the distinct chord lengths, so \(K(A) = (|H| - 1)/2 = k\). The divisibility statement follows from \(H \leq \mathbb{Z}_{m}\). \end{proof}

\begin{cor}[Classification and attainment]\label{res:4.3}
Fix integers \(t \geq 0\) and \(n > \max\{ t,(\binom{t}{2})\}\) with \(n + t \geq 3\) and \(n + t\) odd. Deleting any \(t\) vertices from a regular \((n + t)\)-gon gives a set with maximum reflection defect exactly \(t\) and local maximum \(k = (n + t - 1)/2\). Conversely, every set satisfying Corollary~\ref{res:2.2} is of this form in its minimal regular completion.

\end{cor}

\begin{proof}
Only the forward direction remains. Let \(B\) be the deleted set, with \(|B| = t\). Since \(n > t\), the remaining set \(A\) has density greater than \(1/2\) in \(\mathbb{Z}_{n + t}\); hence \(A - A = \mathbb{Z}_{n + t}\), by intersecting \(A\) and \(A + h\) for each \(h\). Every reflection defect is at most \(t\) by (5). At most \((\binom{t}{2})\) centers can be midpoints of two distinct holes. A diagonal pair of holes cannot have its midpoint in \(A\). Since \(n > (\binom{t}{2})\), some \(a \in A\) has \(r_{B + B}(2a) = 0\), so \(\sigma_{A}(a) = t\). Equation~(3) completes the proof. \end{proof}

Thus the coefficient of \(t\) in the exact hole count is attained for arbitrary patterns of holes. This assertion does not identify the optimal threshold \(Q(t)\).

\section{Sum-cover obstructions and a hierarchy of distance gaps}\label{sum-cover-obstructions-and-a-hierarchy-of-distance-gaps}

In this section \(A\) is translated so that its minimal affine completion is the subgroup \(H\) itself. Put \(B = H\backslash A\) and

\[\rho(B;A) = \min_{a \in A}r_{B + B}(2a).\]

A subset \(B \subseteq H\) is a \emph{sum cover} if \(B + B = H\); repetition in the sum is allowed. We distinguish the sumset \(B + B\) from the dilate \(2B = \{ 2b:b \in B\}\).

\begin{thm}[Exact obstruction and gap identity]\label{res:5.1}
For every nonempty \(A\) in an odd-order abelian group, with the above translation and notation,

\begin{equation}
u - t = \rho(B;A),\quad\quad u = t\  \Leftrightarrow \ B + B \neq H.\label{eq:9}
\end{equation}

The integer \(u - t\) is nonnegative and even. If \(A\) is a subset of an odd regular polygon and \(A - A = H\), then

\begin{equation}
K(A) - k = \frac{u - t}{2} = \frac{1}{2}\min_{s \in H\backslash 2B}r_{B + B}(s).\label{eq:10}
\end{equation}

In particular, the gap is positive exactly when \(B\) is a sum cover. The assumption \(A - A = H\) is automatic when \(n > 4t\), or when \(2n > |H|\).

\end{thm}

\begin{proof}
Taking the maximum in (5) proves \(u - t = \rho(B;A)\). Every element of \(2B\) is in \(B + B\) by a diagonal representation. Since doubling is a bijection, \(2A = H\backslash 2B\). Therefore \(B + B \neq H\) precisely when some \(2a\), \(a \in A\), has no representation, which is equivalent to \(\rho(B;A) = 0\). The parity assertion was proved in Lemma~\ref{res:4.1}. When \(A - A = H\), all chord classes of the completion occur and \(K(A) = (n + u - 1)/2\), whereas Lemma~\ref{res:3.1} gives \(k = (n + t - 1)/2\). Their difference is (10). The final assertions follow from Lemma~\ref{res:3.2} and, respectively, the nonempty intersection of \(A\) and \(A + h\) for every \(h \in H\). \end{proof}

\begin{prop}[Representation budget]\label{res:5.2}
With \(\Delta = u - t\), one has

\begin{equation}
\Delta n \leq u(u - 1).\label{eq:11}
\end{equation}

If \(B\) is a sum cover and

\begin{equation}
4|H| > u^{2} + 3u,\label{eq:12}
\end{equation}

then \(\Delta = 2\). If in addition \(A - A = H\) in the polygon setting, then \(K(A) = k + 1\).

\end{prop}

\begin{proof}
The \(u^{2}\) ordered sums from \(B\) are distributed among \(2A\) and \(2B\). Each of the \(n\) elements of \(2A\) has at least \(\Delta\) representations; the \(u\) distinct elements of \(2B\) each have at least their diagonal representation. Hence \(u^{2} \geq \Delta n + u\). If \(B\) is a sum cover, \(\Delta\) is a positive even integer. If \(\Delta \geq 4\), the same inequality gives \(u^{2} \geq 4(|H| - u) + u = 4|H| - 3u\), contradicting (12). Thus \(\Delta = 2\), and the last claim follows from (10). \end{proof}

\begin{thm}[Finite hierarchy]\label{res:5.3}
Let \(j \geq 1\) be an integer. If

\begin{equation}
n > \max\{ 6t,\frac{(t + 2j)(t + 2j - 1)}{2j}\},\label{eq:13}
\end{equation}

then \(A - A = H\) and \(u \leq t + 2j - 2\). In the polygon setting, \(K(A) - k \leq j - 1\).

\end{thm}

\begin{proof}
Proposition~\ref{res:4.2} gives \(A - A = H\) and \(u < n/2\). Suppose \(u \geq t + 2j\). On the interval \(\lbrack 0,n/2\rbrack\), the function \(f(v) = v(n - v + 1)\) is increasing. Consequently,

\[nt \geq f(u) \geq f(t + 2j) = (t + 2j)(n - t - 2j + 1).\]

Rearrangement gives \(2jn \leq (t + 2j)(t + 2j - 1)\), contrary to (13). Thus \(u < t + 2j\), and parity gives the claimed bound. Equation~(10) gives its geometric interpretation. \end{proof}

\begin{cor}[The one-distance-gap range]\label{res:5.4}
For an odd regular polygon, if

\[n > \max\{ 6t,\frac{(t + 4)(t + 3)}{4}\},\]

then \(K(A) \in \{ k,k + 1\}\). After translating the minimal completion to \(H\), the first case holds if and only if \(B + B \neq H\); the second holds if and only if \(B + B = H\).

\end{cor}

\begin{proof}
Apply Theorem~\ref{res:5.3} with \(j = 2\) and then Theorem~\ref{res:5.1}. \end{proof}

This is an exact dichotomy with an explicit arithmetic certificate, not a classification of all minimum sum covers. The case \(j = 1\) recovers Theorem~\ref{res:2.1}. No optimality is asserted for the thresholds in the hierarchy.

\begin{rem}[Counts and equality cases]\label{res:5.5}
Under Theorem~\ref{res:2.1}, at most \((\binom{t}{2})\) points of \(A\) are midpoints of two distinct holes. Thus \(|\{ a \in A:\sigma_{A}(a) = t\}| \geq n - (\binom{t}{2})\). This is attainable: take holes \(1,4,\ldots,4^{t - 1}\) in a sufficiently large odd cyclic group. Their unordered pair sums are distinct and no off-diagonal sum is twice a hole. The case \(t = 0\) is immediate. In the polygon setting, inclusion--exclusion gives \(r_{A}(h) = n - t + r_{B}(h) \geq n - t\) for \(h \in H\backslash\{ 0\}\). Each chord class therefore occurs in at least \(n - t\) unordered vertex pairs.

The equality \(nt = u(n - u + 1)\) holds precisely when every defect equals \(t\) and \(T_{3}(B) = u\), since its slack is \(\sum_{a \in A}^{}{(t - \sigma_{A}(a))} + T_{3}(B) - u\).

\end{rem}

\begin{rem}[Missing difference classes]\label{res:5.6}
For every subset of an odd regular polygon, without assuming \(A - A = H\),

\[K(A) - k = \frac{u - t - (|H| - |A - A|)}{2}.\]

Indeed \(K(A) = (|A - A| - 1)/2\) and \(k = (n + t - 1)/2\). Thus a sum-cover obstruction alone need not force a positive distance gap.

\end{rem}

\begin{rem}[A prime-order example]\label{res:5.7}
In \(H = \mathbb{Z}_{13}\) take \(B = \{ 0,1,2,6,9\}\) and \(A = H\backslash B\). The sums \(0 + 0,0 + 1,0 + 2,1 + 2,2 + 2,9 + 9,0 + 6,1 + 6,2 + 6,0 + 9,1 + 9,2 + 9,6 + 6\) give residues \(0,1,\ldots,12\). Hence \(B + B = H\). Here \(n = 8\), \(u = 5\), \(A - A = H\) by density, and \(4|H| = 52 > u^{2} + 3u = 40\). Proposition~\ref{res:5.2} gives \(t = 3\), so \(k = 5\), \(K(A) = 6\). Prime order therefore does not exclude completion failure.

\end{rem}

\section{The extremal threshold}\label{the-extremal-threshold}

For the remaining sections, a cyclic set is identified with a subset of its minimal completion \(\mathbb{Z}_{M}\), after translation. Thus \(\langle A - A\rangle = \mathbb{Z}_{M}\). Define

\begin{equation}
F(t) = \max(\{ 0\} \cup \{\,|A|:|A| \geq 2,\ M\text{ odd},\ t(A) = t,\ M - |A| > t\,\}).\label{eq:14}
\end{equation}

Theorem~\ref{res:2.1} ensures finiteness. Consequently \(n > F(t)\) is the optimal strict integer threshold for the conclusion \(u = t\), with the defect fixed \emph{exactly} at \(t\). In the cyclic setting this conclusion also implies \(K = k\): indeed \(K \leq (M - 1)/2 = (n + t - 1)/2 = k\), and always \(K \geq k\). In particular all differences occur. A failure of \(u = t\) need not itself force \(K > k\) when \(A - A \neq \mathbb{Z}_{M}\), as a two-point set illustrates.

\section{The asymptotic completion coefficient}\label{the-asymptotic-completion-coefficient}

\subsection{A cyclic family of obstructions}\label{a-cyclic-family-of-obstructions}

The following construction supplies a positive lower bound and nonempty extremal classes at arbitrarily large scales. The coefficient \(1/3\) is used to establish the identification below and is not claimed to be the best lower bound.

\begin{thm}[A three-fibre odd cyclic obstruction]\label{res:7.1}
Let \(q \geq 7\) with \(q \equiv 1\ (mod\ 6)\), put \(r = q - 2\), and set

\[H = \mathbb{Z}_{q} \times \mathbb{Z}_{r} \times \mathbb{Z}_{3},\quad\quad M = 3q(q - 2).\]

Define the disjoint fibres

\[\begin{aligned}
B_{0} & = \{(x,0,0):0 \leq x < q\}, \\
B_{1} & = \{(0,y,1):0 \leq y < r\}, \\
B_{2} & = \{(i\ mod\ q,i\ mod\ r,2):0 \leq i < q\},
\end{aligned}\]

and let \(B = B_{0} \cup B_{1} \cup B_{2}\), \(A = H\backslash B\). The group \(H\) is cyclic of odd order and \(B + B = H\). Identifying it with the vertex labels of a regular \(M\)-gon gives

\begin{equation}
\begin{aligned}
u & = 3q - 2, & t & = 3q - 4, & n & = 3q^{2} - 9q + 2, \\
k & = (M - 3)/2, & K(A) & = (M - 1)/2, & K(A) - k & = 1.
\end{aligned}
\label{eq:15}
\end{equation}

In particular,

\begin{equation}
\frac{n}{t^{2}} \rightarrow \frac{1}{3}.\label{eq:16}
\end{equation}

No universal sufficient leading coefficient for exact completion can be smaller than \(1/3\).

\end{thm}

\begin{proof}
The three moduli \(q,r,3\) are pairwise coprime: \(q,r\) are odd and differ by \(2\), while \(q \equiv 1\) and \(r \equiv 2\ (mod\ 3)\). The Chinese remainder theorem therefore identifies \(H\) with \(\mathbb{Z}_{M}\). Distinct third coordinates separate the fibres, and the first coordinate separates the \(q\) elements of \(B_{2}\). Thus \(u = q + r + q = 3q - 2\).

Every element \((x,y,z)\) is a mixed sum of two fibres. If \(z = 0\), choose the unique \(i \in \{ 0,\ldots,q - 1\}\) with \(i = x\ (mod\ q)\); then

\[(x,y,0) = (0,y - i,1) + (i,i,2) \in B_{1} + B_{2}.\]

If \(z = 1\), write \((x,y,1) = (x,0,0) + (0,y,1)\). If \(z = 2\), choose \(i \in \{ 0,\ldots,r - 1\}\) with \(i = y\ (mod\ r)\) and write

\[(x,y,2) = (x - i,0,0) + (i,i,2) \in B_{0} + B_{2}.\]

Coordinates in these expressions are reduced in their respective groups. This proves \(B + B = H\).

We have \(M > 2u\) for \(q \geq 7\), so \(A - A = H\) by density. Moreover,

\[4M - u^{2} - 3u = 3q^{2} - 21q + 2 > 0\]

for every \(q \geq 7\). Proposition~\ref{res:5.2} gives \(u - t = 2\). Subtraction yields (15), and (10) yields the distance counts. Equation~(16) follows directly. For every \(c < 1/3\) and every function \(g(t) = o(t^{2})\), these counterexamples eventually satisfy \(n > ct^{2} + g(t)\); such a condition cannot imply \(u = t\). \end{proof}

\subsection{Identification of the asymptotic threshold}\label{identification-of-the-asymptotic-threshold}

We now isolate precisely what is and is not determined about the leading constant. Define

\begin{align}
\alpha_{\mathrm{poly}} & = \lim_{T \rightarrow \infty}\sup\{\frac{|A|}{t{(A)}^{2}}:\begin{gathered}
A \subseteq \mathbb{Z}_{m},\ m \geq 3\text{ odd},\ |A| \geq 2, \\
t(A) \geq T,\ u(A) > t(A)
\end{gathered}\},\label{eq:17} \\
\beta_{\mathrm{odd}} & = \lim_{U \rightarrow \infty}\sup\{\frac{M}{{|B|}^{2}}:\begin{gathered}
M\text{ odd},\ B \subseteq \mathbb{Z}_{M}, \\
|B| \geq U,\ B + B = \mathbb{Z}_{M}
\end{gathered}\}.\label{eq:18}
\end{align}

Here \(u(A)\) is always measured in the minimal affine completion of \(A\). Both limits exist as decreasing limits of suprema; Theorem~\ref{res:7.1} makes their defining classes nonempty for arbitrarily large cutoffs. Theorem~\ref{res:2.1} bounds \(\alpha_{\mathrm{poly}}\) by \(1/2\), while counting unordered pairs with repetition gives

\begin{equation}
M \leq \frac{u(u + 1)}{2}\quad\quad(u = |B|,\ B + B = \mathbb{Z}_{M}),\label{eq:19}
\end{equation}

so \(\beta_{\mathrm{odd}} \leq 1/2\).

\begin{thm}[Equality of the extremal constants]\label{res:7.2}
The optimal asymptotic completion coefficient is exactly the odd cyclic sum-cover efficiency:

\begin{equation}
\alpha_{\mathrm{poly}} = \beta_{\mathrm{odd}},\quad\quad\frac{1}{3} \leq \alpha_{\mathrm{poly}} \leq \frac{1}{2}.\label{eq:20}
\end{equation}

More explicitly, for every \(\varepsilon > 0\) and all sufficiently large \(t\), the condition \(n > (\alpha_{\mathrm{poly}} + \varepsilon)t^{2}\) implies \(u = t\), whereas there are counterexamples with arbitrarily large \(t\) and \(n > (\alpha_{\mathrm{poly}} - \varepsilon)t^{2}\).

\end{thm}

\begin{proof}
The three-fibre family proves that both constants are at least \(1/3\). Choose a sequence of counterexamples with \(t \rightarrow \infty\) and \(n/t^{2} \rightarrow \alpha_{\mathrm{poly}}\). Since the limit is positive, \(n > 6t\) eventually. Proposition~\ref{res:4.2}, or the preceding hull bound, gives

\[u \leq \frac{2nt}{n - 2t} = O(t).\]

By (11),

\[0 < u - t \leq \frac{u(u - 1)}{n} = O(1).\]

Thus \(u/t \rightarrow 1\) and \(M = n + u\) satisfies

\[\frac{M}{u^{2}} - \frac{n}{t^{2}} \rightarrow 0.\]

After translation and a cyclic identification, Theorem~\ref{res:5.1} supplies a sum cover of \(\mathbb{Z}_{M}\). Since \(u \rightarrow \infty\), this sequence shows \(\alpha_{\mathrm{poly}} \leq \beta_{\mathrm{odd}}\).

Conversely, choose sum covers with \(u \rightarrow \infty\) and \(M/u^{2} \rightarrow \beta_{\mathrm{odd}} \geq 1/3\). For \(A = \mathbb{Z}_{M}\backslash B\), we have \(u/M \rightarrow 0\), so \(A\) is eventually nonempty, has more than \(M/2\) elements, and satisfies \(A - A = \mathbb{Z}_{M}\). In particular its minimal completion is the whole group. Since the limit exceeds \(1/4\), eventually \(4M > u^{2} + 3u\). Proposition~\ref{res:5.2} then gives \(t = u - 2\). It follows that

\[\frac{n}{t^{2}} = \frac{M - u}{{(u - 2)}^{2}} \rightarrow \beta_{\mathrm{odd}}.\]

These are polygon counterexamples, proving the reverse inequality. The quantitative interpretation with \(\varepsilon\) is precisely the defining decreasing supremum in (17). \end{proof}

Theorem~\ref{res:7.2} and the unordered-pair bound show that \(\alpha_{\mathrm{poly}} = 1/2\) if and only if there are odd \(M \rightarrow \infty\) and sum covers with \(|B| = (\sqrt{2} + o(1))\sqrt{M}\).

This is the odd-modulus version of Green's Problem~33 \cite{green}. A positive answer here would answer that unrestricted-modulus question positively; a negative answer for odd moduli alone would not settle its even-modulus cases. No numerical evaluation of \(\beta_{\mathrm{odd}}\) is proved here.

\begin{rem}[Two consequences of the constructions]\label{res:7.3}
Substituting \(q = (t + 4)/3\) in Theorem~\ref{res:7.1} yields \(F(t) \geq (t^{2} - t - 14)/3\) for \(t \equiv 17\ (mod\ 18)\), \(t \geq 17\). This is a lower bound from that construction, not an optimality assertion.

For any sequence of completion counterexamples with \(t \rightarrow \infty\) and \(n/t^{2} \rightarrow \alpha_{\mathrm{poly}}\), the preceding proof gives \(u = t + O(1)\) and \(M/u^{2} \rightarrow \alpha_{\mathrm{poly}} \geq 1/3\). Thus eventually \(4M > u^{2} + 3u\) and \(n > 4t\). Proposition~\ref{res:5.2} and Theorem~\ref{res:5.1} imply \(u - t = 2\) and \(K(A) - k = 1\).

\end{rem}

\section{Scope}\label{scope}

The completion theorem applies to finite abelian groups of odd order. Its geometric interpretation concerns subsets of an odd regular polygon on a common circle. The observed set need not inherit the symmetries of its completion, since arbitrary patterns of \(t\) missing vertices occur in Corollary~\ref{res:4.3}.

Odd order is essential to the argument. For even ambient order the antipodal fixed point changes the pinned-distance identity to

\begin{equation}
2d_{A}(a) = n - 1 + \sigma_{A}(a) + \eta_{A}(a),\quad\quad\eta_{A}(a) = \mathbf{1}_{\{ a + m/2 \in A\}}.\label{eq:21}
\end{equation}

Indeed, besides the fixed centre, its antipode contributes one distance and no reflection defect; all other reflection orbits have size two. No even-order completion theorem is asserted. The exact value of \(\alpha_{\mathrm{poly}}\) and the existence of a limit \(F(t)/t^{2}\) remain open. The numerical inequality in Theorem~\ref{res:7.2} is illustrative and is not a claim to the best available bound.

The analogous asymptotic identification for all finite abelian groups of odd order is proved in Section~S10. The remaining questions include numerical thresholds beyond the retained small cases, existence of \(\lim F(t)/t^{2}\), and stronger bounds or efficient infinite obstruction families under prime-order completion. The example above settles only existence of a prime-order obstruction.

\section*{Acknowledgement of AI assistance}
An AI tool was used minimally and strictly for language polishing and grammar checking; no mathematical content, proofs, or ideas were generated by AI.
\clearpage
\section*{Supplementary material}
\markboth{J. LEE}{SUPPLEMENT: COMPLETION IN ODD REGULAR POLYGONS}
This Supplement contains additional results, finite classifications, and computational records. The main article contains the proofs of its principal theorems. Numbers prefixed by S refer to this Supplement; other numbered references are to the main article. A common bibliography follows the Supplement. The computational source files, certificates, and recorded outputs are supplied in the ancillary directory \texttt{anc/code/}.
\setcounter{section}{0}
\setcounter{equation}{0}
\setcounter{table}{0}
\setcounter{figure}{0}
\renewcommand{\thesection}{S\arabic{section}}
\renewcommand{\theequation}{S\arabic{equation}}
\renewcommand{\thetable}{S\arabic{table}}
\renewcommand{\thefigure}{S\arabic{figure}}
\renewcommand{\theHsection}{sup.\arabic{section}}
\renewcommand{\theHequation}{sup.\arabic{equation}}
\renewcommand{\theHtable}{sup.\arabic{table}}
\renewcommand{\theHfigure}{sup.\arabic{figure}}
\section{Finite thresholds: summary}\label{sup:sec:finite_summary}
\begin{suptheorem}[Finite computability and the retained small cases]\label{sup:thm:finiteoverview}
For every nonnegative integer $t$, the value $F(t)$ is computable by a terminating finite algorithm. If $F(t)>2t$, all extremizers can also be listed up to affine cyclic equivalence. The following values hold:
\[
 (F(0),F(1),F(2),F(3),F(4),F(5))=(0,2,5,8,13,14).
\]
The extremizers are completely classified for $t\leq4$. No complete extremal classification is asserted for $t=5$.
\end{suptheorem}
\begin{proof}
The analytic finite reduction is proved in Supplement Theorems~\ref{sup:thm:finiteM} and~\ref{sup:thm:computable}. The case $t=0$ follows from Theorem~\ref{res:2.1}. The remaining values and classifications follow from Supplement Theorems~\ref{sup:thm:F1}, \ref{sup:thm:F2}, \ref{sup:thm:F34} and~\ref{sup:thm:F5}. From $t=2$ onward these are computer-assisted proofs, with globally justified finite search bounds; they are not conclusions from a bounded sample of polygons.
\end{proof}

\begin{table}[ht]
\centering
\begin{tabular}{@{}crll@{}}\toprule
$t$ & $F(t)$ & Extremal classification & Proof component\\\midrule
0&0&No counterexamples&Analytic\\
1&2&All odd $M\geq5$; one affine class per $M$&Analytic\\
2&5&Two affine classes; four dihedral classes&Computer-assisted\\
3&8&One affine class; three dihedral classes&Computer-assisted\\
4&13&One affine class; three dihedral classes&Computer-assisted\\
5&14&One witness; full classification not asserted&Computer-assisted\\\bottomrule
\end{tabular}
\caption{The retained finite results. Affine equivalence is $A\mapsto rA+s$ with $r$ a unit; dihedral equivalence allows only $r=\pm1$.}
\label{sup:tab:finiteoverview}
\end{table}

The structural endpoint of the finite classification is $t=4$. The value at $t=5$ is retained as a short additional entry, with its complete proof confined to the Supplement. These computations illustrate the theory but are not used to prove the general asymptotic results. Computability does not supply a practical algorithm or a closed numerical formula.

\section{Finite reduction and defects one and two}\label{sup:sec:small}
Use the notation and definition of $F(t)$ from the main manuscript. Every cyclic set below is identified with its minimal completion after translation.

\subsection{Bounding the completion without bounding the ambient polygon}
\begin{suptheorem}[Finite completion bound]\label{sup:thm:finiteM}
Let $A$ have $n\geq2$ elements in a finite abelian group of odd order, with maximum reflection defect $t$. If $n>2t$, its minimal completion has order
\begin{equation}\label{sup:eq:finiteM}
 M\leq \left\lfloor6^{(n-1)/2}\right\rfloor.
\end{equation}
Moreover, for $D=A-A$,
\begin{equation}\label{sup:eq:finiteD}
 2D=D=-D,\qquad
 |D|\leq 1+\left\lfloor\frac{n(n-1)}{n-2t}\right\rfloor.
\end{equation}
\end{suptheorem}
\begin{proof}
Put $p=(n-1-t)/2$. By Lemma~\ref{res:3.1}, this is an integer. The hypotheses imply $p\geq1$. Label the points $a_0,\ldots,a_{n-1}$. At each $a_i$ choose $p$ disjoint reflection pairs $\{a_b,a_c\}$, so that $2a_i=a_b+a_c$. Form the integer relation rows $2e_i-e_b-e_c$ in $\Z^n$.

Construct a directed graph by joining $i$ to all $2p=n-1-t$ chosen endpoints. A closed strongly connected component has at least $2p+1=n-t$ vertices. Since $n>2t$, there can be only one such component; every vertex reaches it. If a real vector satisfies all the chosen relations, then its value at each vertex is the average of its values at the outgoing neighbors. On a closed strongly connected component the maximum principle makes the vector constant. A global maximum and a global minimum both propagate along directed paths to that same component. They therefore coincide, and the vector is constant everywhere. Thus the relation matrix has rank $n-1$ over $\mathbb Q$.

Translate $a_0$ to zero and delete column $0$. Choose $n-1$ independent rows, forming a nonsingular integer matrix $R$. The assignment of its coordinate generators to $a_1,\ldots,a_{n-1}$ induces a surjection
\[
 \Z^{n-1}/\langle\text{rows of }R\rangle\ \longrightarrow\ H.
\]
Its domain has order $|\det R|$, so $M$ divides this determinant. Each row has Euclidean norm at most $\sqrt6$. Hadamard's determinant inequality proves \eqref{sup:eq:finiteM}.

The first part of the proof of Lemma~\ref{res:3.2} uses only $n>2t$ to give $2D=D$ and $r_A(h)\geq n-2t$ for $h\in D$. Since $r_A(0)=n$,
\[
 n(n-1)=\sum_{h\in D\setminus\{0\}}r_A(h)\geq (|D|-1)(n-2t),
\]
which proves \eqref{sup:eq:finiteD}. Symmetry $D=-D$ is automatic.
\end{proof}

The strict inequality $n>2t$ cannot simply be removed: the set $\{0,1\}\subseteq\Z_M$ has $n=2$, $t=1$, and unbounded minimal completion order as odd $M$ varies. The bound \eqref{sup:eq:finiteM} is deliberately coarse; its purpose is to turn a search over all odd polygons into a certified finite search.

\begin{suptheorem}[Effective determination at fixed defect]\label{sup:thm:computable}
There is a terminating finite algorithm that computes $F(t)$ for every integer $t\geq0$. If $F(t)>2t$, it also lists all extremal sets up to affine cyclic equivalence. This assertion is about computability, not a practical complexity bound or a closed formula for $F(t)$.
\end{suptheorem}
\begin{proof}
For $t=0$, Theorem~\ref{res:2.1} applies to every admissible $n\geq2$, so $F(0)=0$. For odd $t\geq1$, take in $\Z_{5t}$ the two full residue classes $0,1$ modulo $5$. They give $n=2t$, actual maximum defect $t$, and $u=3t>t$. For even $t\geq4$, put $h=t-1$ and take
\[
 A=\{x\in\Z_{5h}:x\equiv0\text{ or }1\pmod5\}\cup\{3\}.
\]
Then $n=2t-1$. At every old point the defect is $h+1=t$: reflection preserves its own residue class, misses the other class, and sends the added point outside $A$. At the new point $3$, reflection interchanges the two full classes and fixes $3$, giving defect zero. Also $u=3h-1>h+1=t$. These sets contain $0,1$, so their minimal completion is the whole group. For $t=2$, the directly checkable set $\{0,1,2,3,6\}\subseteq\Z_9$ has $n=5$, $t=2$, and $u=4$.

For odd $t$ the largest possible size at most $2t$ is already realized by a counterexample. For even $t\geq4$, parity $n\equiv t+1\pmod2$ makes $2t-1$ the largest admissible size at most $2t$, and the displayed example realizes it. The $t=2$ example already exceeds $2t$. Thus smaller sizes cannot affect the maximum once these witnesses are recorded. It remains to test the finitely many integers
\[
 2t<n\leq Q(t),\qquad n\equiv t+1\pmod2.
\]
For each, enumerate all odd $M$ between $n$ and the bound in \eqref{sup:eq:finiteM}, and all $n$-element subsets of $\Z_M$ containing $0$. Retain precisely those generating the full completion, having actual defect $t$, and satisfying $M-n>t$. The maximum of the retained sizes and the recorded witness size is $F(t)$. The same list, quotiented by $A\mapsto rA+s$ with $\gcd(r,M)=1$, gives the asserted extremal classification when $F(t)>2t$.
\end{proof}

\subsection{The first two positive defects}
\begin{suptheorem}[Defect one]\label{sup:thm:F1}
One has $F(1)=2$. Every counterexample has two points and minimal completion order an odd integer $M\geq5$. For each such $M$ there is one affine class, represented by $\{0,1\}\subseteq\Z_M$. Up to rotation and reflection, the representatives are $\{0,d\}$ with $1\leq d\leq(M-1)/2$ and $\gcd(d,M)=1$.
\end{suptheorem}
\begin{proof}
Parity makes $n$ even. If $n\geq8$, Theorem~\ref{res:2.1} applies because $Q(1)=6$. For $n=6$, Lemma~\ref{res:3.2} gives $M\leq9$, so $M=7$ or $9$. The latter would give $u=3$ and contradict the representation budget $(u-t)n\leq u(u-1)$ from Proposition~\ref{res:5.2}. Hence $u=1$.

For $n=4$, choose a reflection pair at one point and translate that point to $0$, writing $A=\{0,b,-b,c\}$. If $3b=0$, the reflection pair required at $c$ gives $2c\in\langle b\rangle$. Since the quotient has odd order, this forces $c\in\langle b\rangle$, a contradiction to four distinct points. If $3b\ne0$, the reflection pairs at $b$ and $-b$ require respectively
\[
 c\in\{2b,3b\},\qquad c\in\{-2b,-3b\}.
\]
Odd order rules out $4b=0$, and $6b=0$ would imply $3b=0$. Thus $5b=0$, with $c=2b$ or $3b$. The minimal completion has order $5$, again giving $u=1$. Finally, every two-point set has defect $1$. After identifying its minimal cyclic completion, it is affine equivalent to $\{0,1\}$; it fails exactness precisely when $M-2>1$.
\end{proof}

The next lemma is explicitly computer-assisted. Its finite reduction is proved here, and two different exact enumeration programs and their actual output are supplied. No arbitrary cutoff on polygon order is assumed.

\begin{suplemma}[Finite classification at five and seven points]\label{sup:lem:smallclass}
For $n\in\{5,7\}$, the sets of $n$ points in odd cyclic groups with maximum defect at most $2$ are, in their minimal completions, exactly the affine classes in Table~\ref{sup:tab:smallclasses}.
\end{suplemma}
\begin{table}[ht]
\centering
\begin{tabular}{@{}rrlrr@{}}\toprule
$n$ & $M$ & Representative $A$ & $t(A)$ & $u$\\\midrule
5&5&$\{0,1,2,3,4\}$&0&0\\
5&7&$\{0,1,2,3,4\}$&2&2\\
5&9&$\{0,1,2,3,6\}$&2&4\\
5&11&$\{0,1,2,4,7\}$&2&6\\\midrule
7&7&$\{0,1,2,3,4,5,6\}$&0&0\\
7&9&$\{0,1,2,3,4,5,6\}$&2&2\\
7&9&$\{0,1,2,3,4,6,7\}$&2&2\\\bottomrule
\end{tabular}
\caption{Complete affine classification with maximum defect at most two.}
\label{sup:tab:smallclasses}
\end{table}
\begin{proof}[Computer-assisted proof]
Theorem~\ref{sup:thm:finiteM} bounds $M$ by $36$ for $n=5$ and by $216$ for $n=7$. Equation~\eqref{sup:eq:finiteD}, with the upper bound $t\leq2$, bounds $|D|$ by $21$ and $15$, respectively. In each odd $\Z_M$, the set $D$ contains $0$ and is a union of orbits under $x\mapsto2x$ and $x\mapsto-x$. Enumerate all such unions with the stated size bound and $\langle D\rangle=\Z_M$. After translating a point to zero, one has $A\subseteq D$. Enumerate all $n$-element subsets of each candidate $D$ containing $0$, test every reflection defect, and retain a subset only when its actual difference set equals $D$. Finally identify sets under unit multiplication and translation. Every eligible set enters this enumeration, and every retained set is directly checked.

The standard-library program \texttt{verify\_small\_defects.py} implements precisely these steps using integers. For $n=5$ it tests $38$ candidate difference sets and $61{,}653$ candidate subsets; for $n=7$ the counts are $81$ and $144{,}980$. Its complete output is Table~\ref{sup:tab:smallclasses}, recorded in \texttt{small\_defect\_verification.json}.

An independent enumeration in \texttt{classify\_relations.py} chooses the reflection pairs at each center, normalizes the pairs at $0$, and forms their integer relation lattice. Column Hermite normal form removes duplicate lattices. Whenever the lattice has full rank, all solutions of its transpose congruences in $(\mathbb Q/\Z)^d$ are enumerated using exact rational arithmetic, and all distinct-point, odd-order and defect conditions are tested directly. Full rank for a complete pattern follows from the maximum-principle argument in Theorem~\ref{sup:thm:finiteM}. The runs visit $255$ and $2{,}307$ lattice states for $n=5,7$, respectively, and give exactly the same affine classes. This second computation is a cross-check; the first enumeration already suffices with the proved finite bounds.
\end{proof}

\begin{suptheorem}[Defect two and all extremal polygons]\label{sup:thm:F2}
One has $F(2)=5$. At the threshold there are exactly two affine classes of counterexamples:
\[
 \{0,1,2,3,6\}\subseteq\Z_9,\qquad
 \{0,1,2,4,7\}\subseteq\Z_{11}.
\]
Up to rotations and reflections of their minimal regular polygons, these split into the four classes in Table~\ref{sup:tab:dihedral}. All sets of actual defect two with more than five vertices have $u=2$ and $K=k$.
\end{suptheorem}
\begin{proof}
Parity makes $n$ odd. For $n\geq13$, Theorem~\ref{res:2.1} applies, since $Q(2)=12$. If $n=9$ or $11$, Lemma~\ref{res:3.2} gives $M\leq n^2/(n-4)$, and therefore $u\leq6$. By parity and $u\geq t$, the possibilities are $u=2,4,6$. The representation budget excludes $u=4$, since $2n>12$, and excludes $u=6$, since $4n>30$. Thus $u=2$. For $n=7$, Lemma~\ref{sup:lem:smallclass} gives $u=2$ whenever the actual defect is two. Its two five-point counterexamples prove the lower bound and exhaust the extremal affine classes. Enumerating the affine images of those representatives and quotienting by $A\mapsto\pm A+s$ gives Table~\ref{sup:tab:dihedral}; the affine orbits have $27$ and $22$ members, splitting into three and one dihedral orbits, respectively. These orbit calculations are included in the independent verifier.
\end{proof}

\begin{table}[ht]
\centering
\begin{tabular}{@{}rlrrrr@{}}\toprule
$M$ & Representative $A$ & $t$ & $u$ & $k$ & $K$\\\midrule
9&$\{0,1,2,3,6\}$&2&4&3&4\\
9&$\{0,1,2,4,7\}$&2&4&3&4\\
9&$\{0,1,3,5,6\}$&2&4&3&4\\
11&$\{0,1,2,4,7\}$&2&6&3&5\\\bottomrule
\end{tabular}
\caption{All extremal defect-two configurations up to rotation and reflection. Affine unit multiplication preserves distance counts but need not preserve individual chord lengths or Euclidean congruence.}
\label{sup:tab:dihedral}
\end{table}

The $11$-point-completion example is affine equivalent to the classical cyclic $(11,5,2)$ Paley difference set: every nonzero difference has two ordered representations. Such difference sets are established objects; see, for example, \cite[Table~1]{custic}. We claim neither a new construction of that set nor priority for the small configurations themselves. The result established here is their exhaustive role in this fixed-defect completion threshold, with the computer-assisted component disclosed.

\subsection{Lifting the small extremizers}
\begin{supproposition}[Full-fibre lifting]\label{sup:prop:lift}
Let $A_0\subseteq\Z_M$ have minimal completion $\Z_M$, with parameters $(n_0,t_0,u_0)$. For an odd positive integer $h$, let $\pi:\Z_{hM}\to\Z_M$ be reduction modulo $M$ and put $A=\pi^{-1}(A_0)$. Then $A$ has minimal completion $\Z_{hM}$ and
\[
 n=hn_0,\qquad t=ht_0,\qquad u=hu_0.
\]
In particular, for every odd $h\geq1$,
\begin{equation}\label{sup:eq:liftlower}
 F(2h)\geq5h.
\end{equation}
The two seeds of Theorem~\ref{sup:thm:F2} give, respectively, $K-k=h$ and $K-k=2h$.
\end{supproposition}
\begin{proof}
Reflection commutes with $\pi$, so each missing reflected residue contributes exactly $h$ missing vertices. Hence $\sigma_A(a)=h\sigma_{A_0}(\pi(a))$. Moreover $A-A$ contains the kernel of $\pi$ and projects onto a generating subset of $\Z_M$. Thus it generates all of $\Z_{hM}$, proving the three parameter identities. Failed completion is preserved. Both five-point seeds have full difference set (by direct enumeration); their lifts do as well, since $A-A=\pi^{-1}(A_0-A_0)$. The gap formula in Theorem~\ref{res:5.1} gives the last assertion.
\end{proof}

The lifted families need not be extremal for $h>1$; they supply structured lower bounds and growing local-to-global gaps.
\section{Complete classifications at defects three and four}\label{sup:sec:extended}
\subsection{An order filter for difference-set enumeration}
\begin{suplemma}[Orders of doubling-negation orbits]\label{sup:lem:orderfilter}
Let $x$ have odd additive order $d>1$ in a cyclic group. Its orbit under doubling and negation has size
\[
 c(d)=2\min\{j\geq1:2^j\equiv1\text{ or }-1\pmod d\}.
\]
Suppose $D=-D=2D$, $0\in D$, $|D|\leq L+1$, and $\langle D\rangle=\Z_M$. Then $M$ is the least common multiple of some distinct odd integers $d_1,\ldots,d_s>1$ satisfying
\[
 \sum_{i=1}^s c(d_i)\leq L,\qquad
 d_i\mid2^{j_i}-1\text{ or }2^{j_i}+1,\quad j_i\leq\lfloor L/2\rfloor.
\]
\end{suplemma}
\begin{proof}
The orbit consists of the residues $\pm2^j x$. Before the first return up to sign, all these signed residues are distinct. If the first return is to $-x$, doubling continues through the negatives; if it is to $x$, the negatives form a disjoint companion orbit. In either case the total size is twice the first return time. The nonzero part of $D$ is a disjoint union of these orbits. Choose one orbit for each distinct additive order that occurs. Their combined size is at most $|D|-1$, and their orders have least common multiple $M$, since $D$ generates the cyclic group. The divisibility statements follow from the return equations.
\end{proof}

For each $(n,T)$ with $n>2T$ and $n-1\equiv T\pmod2$, put
\[
 N_n=\lfloor6^{(n-1)/2}\rfloor,\qquad
 L_{n,T}=\left\lfloor\frac{n(n-1)}{n-2T}\right\rfloor.
\]
Every set with actual defect at most $T$ has $M\leq N_n$ and $|D|\leq L_{n,T}+1$. Lemma~\ref{sup:lem:orderfilter} generates all possible moduli from divisors of finitely many numbers $2^j\pm1$, using least common multiples and the orbit-size budget. For each modulus, enumerate the allowed orbit unions $D$. This avoids scanning every residue in every modulus, while preserving completeness.

After translating a point to zero, an admissible $A\setminus\{0\}$ is an $(n-1)$-clique in the graph on $D\setminus\{0\}$ whose edges join $x,y$ when $x-y\in D$. Enumerate these cliques, test every reflection defect, and require that the actual difference set equal $D$. The final equivalence reduction is by affine unit maps. The programs use only exact integer tests.

\begin{suplemma}[Finite enumeration for defects three and four]\label{sup:lem:extendedclass}
The enumeration above gives the following conclusions.
\begin{enumerate}[leftmargin=*]
\item For $n=8$, $T=3$, there is exactly one affine class with $u>t(A)$, represented by
\[
 A_3=\{0,1,2,3,4,5,8,10\}\subseteq\Z_{13};\quad t(A_3)=3,\ u=5.
\]
For $n=10,12$ and actual defect three, all classes have $u=3$.
\item For $n=13$, $T=4$, there is exactly one affine class with $u>t(A)$, represented by
\[
 A_4=\{0,1,2,3,4,5,6,7,9,12,13,14,17\}\subseteq\Z_{19};\quad t(A_4)=4,\ u=6.
\]
For $n=15$ and actual defect four, all classes have $u=4$.
\item No $17$-element subset of $\Z_{31}$ has maximum defect at most four.
\end{enumerate}
\end{suplemma}
\begin{proof}[Computer-assisted proof]
The finite bounds and Lemma~\ref{sup:lem:orderfilter} prove completeness of the candidate generation. The graph search lists all potentially admissible subsets, and direct membership tests enforce all claimed conditions. The complete affine lists are recorded in the accompanying result files; Table~\ref{sup:tab:extendedchecks} summarizes the runs. The last row only concerns the explicitly specified modulus $31$.

The main search additionally prunes a partial clique if it cannot supply the required number $p=(n-1-T)/2$ of reflection pairs at zero. If it already contains $r$ complete pairs, has $s$ selected singleton partners still available, and has $b$ pairs wholly available, then with $v$ selections left it can gain at most
\[
 \min(s,v)+\min\bigl(b,\lfloor(v-\min(s,v))/2\rfloor\bigr)
\]
additional complete pairs. A branch is discarded only if this upper bound plus $r$ is less than $p$.

For the three defect-four rows, the search was repeated with this pair pruning entirely disabled; the accepted sets and affine classes agreed exactly. For the defect-three rows, an independently implemented Python enumeration using direct modulus scans and clique listing also agreed, including the full candidate difference-set universes. The eight-point result additionally agrees with the relation-lattice enumeration. The executable driver is \texttt{verify\_extended.py}; the integer search is \texttt{classify\_exact.cpp}. These computations are part of the proof, not formal proof-assistant certification.
\end{proof}

\begin{table}[ht]
\centering
\begin{tabular}{@{}rrrrrr@{}}\toprule
$n$ & $T$ & Modulus bound & $|D|$ bound & Candidate $D$ & Affine classes\\\midrule
8&3&529&29&1,214&4\\
10&3&3,174&23&1,328&4\\
12&3&19,047&23&1,344&8\\
13&4&46,656&32&20,935&16\\
15&4&279,936&31&21,803&16\\
17&4&31 only&31&1&0\\\bottomrule
\end{tabular}
\caption{Exhaustive finite enumerations. The bounds $T$ include smaller actual defects. In the fifth column each $D$ contains zero and generates its full cyclic group.}
\label{sup:tab:extendedchecks}
\end{table}

\begin{suptheorem}[Sharp thresholds at defects three and four]\label{sup:thm:F34}
The exact thresholds are
\[
 \boxed{F(3)=8,\qquad F(4)=13.}
\]
For each defect the extremizers form one affine class, represented by $A_3$ or $A_4$ in Lemma~\ref{sup:lem:extendedclass}. Each affine class splits into three classes under rotations and reflections. Their complements are listed in Table~\ref{sup:tab:extendedextreme}. The distance parameters are $(k,K)=(5,6)$ for defect three and $(8,9)$ for defect four.
\end{suptheorem}
\begin{proof}
For $t=3$, parity makes $n$ even and Theorem~\ref{res:2.1} gives exactness for $n>18$. For $n=14,16,18$, Lemma~\ref{res:3.2} bounds $u$ by $9$. If $u>3$, its odd possible values are $5,7,9$. Each violates $(u-3)n\leq u(u-1)$, already at $n=14$. Lemma~\ref{sup:lem:extendedclass} handles $n=10,12$ and supplies the eight-point extremizer. This proves $F(3)=8$ and its affine classification.

For $t=4$, parity makes $n$ odd and Theorem~\ref{res:2.1} gives exactness for $n>24$. For $n=19,21,23$, the hull bound gives $u\leq12$. Every even $u\in\{6,8,10,12\}$ violates $(u-4)n\leq u(u-1)$, already at $n=19$. At $n=17$, the hull bound gives $M\leq32$ and hence $u\leq14$; the budget excludes $u=6,8,10,12$, leaving only $u=14$, $M=31$. This is ruled out by the last part of Lemma~\ref{sup:lem:extendedclass}. Its fifteen-point classification and thirteen-point extremizer finish the proof of $F(4)=13$.

The affine orbits of $A_3,A_4$ have $39$ and $57$ members. Quotienting by translations and sign gives three dihedral classes in each case, as checked directly in \texttt{verify\_extended.py}. Both sets have more than half the points of their completions, so their difference sets are full. The pinned-distance identity and $K=(M-1)/2$ give the stated distance parameters.
\end{proof}

\begin{table}[ht]
\centering
\begin{tabular}{@{}rrrl@{}}\toprule
$t$ & $n$ & $M$ & Deleted vertices $B$ (so $A=\Z_M\setminus B$)\\\midrule
3&8&13&$\{0,1,2,6,9\}$\\
3&8&13&$\{0,1,3,5,6\}$\\
3&8&13&$\{0,1,3,7,11\}$\\\midrule
4&13&19&$\{0,1,3,12,14,15\}$\\
4&13&19&$\{0,1,4,5,11,13\}$\\
4&13&19&$\{0,1,4,6,14,16\}$\\\bottomrule
\end{tabular}
\caption{All extremizers at defects three and four, up to rotations and reflections.}
\label{sup:tab:extendedextreme}
\end{table}

The first cover in each block already appears in Haanp\"a\"a's sum-cover table \cite[Table~3]{haanpaa}; these cover constructions are not new. The added conclusion concerns their exhaustive extremal role for a fixed maximum reflection defect, including the exclusion of larger polygon sets and other completion orders. The established cover constructions must be distinguished from the fixed-defect exclusion and classification statements proved here.

\section{The exact threshold at defect five}\label{sup:sec:F5}
The next finite value is also accessible, but requires stronger pruning in the exhaustive search. The result below determines the threshold; it does not classify every extremal fourteen-point configuration.

\begin{suplemma}[The remaining defect-five cases]\label{sup:lem:F5absence}
There is no set with $n\in\{16,18,20\}$, maximum reflection defect at most five, and more than five holes in its minimal odd cyclic completion. There is also no $22$-element set in $\Z_{39}$ with maximum defect at most five.
\end{suplemma}
\begin{proof}[Computer-assisted proof]
For the first three sizes, Theorem~\ref{sup:thm:finiteM} and Lemma~\ref{sup:lem:orderfilter} give the complete candidate difference sets. Discard moduli with $u=M-n\leq5$ or $(u-5)n>u(u-1)$; these cannot satisfy the claimed conditions. The resulting finite lists have the sizes in Table~\ref{sup:tab:F5checks}. The $22$-point case uses only $D=\Z_{39}$, as justified in the theorem below.

We again enumerate cliques in the graph on $D\setminus\{0\}$. Four additional reductions are safe. First, repeatedly discard a vertex with fewer than $n-2$ neighbors among the remaining vertices; such a vertex cannot lie in an $(n-1)$-clique. Second, any proper coloring of a candidate graph bounds its clique size from above; greedy color classes used here are independent sets. Third, if $D$ contains a unit, an affine map can send the endpoints of a unit difference to $0,1$. It is therefore sufficient to search sets containing $0,1$ among candidate $D$ containing $1$. The candidate universe is invariant under multiplication by units. If $D$ has no unit, the search keeps the normalization at $0$ alone; generating sets without a unit difference are not omitted.

Fourth, apply the reflection-pair upper bound from the proof of Lemma~\ref{sup:lem:extendedclass} at \emph{every selected center}. With $v$ remaining choices, $r$ already complete pairs, $s$ available partners of selected singletons, and $b$ pairs wholly available, the final number of complete pairs is at most
\[
 r+\min(s,v)+\min\bigl(b,\lfloor(v-\min(s,v))/2\rfloor\bigr).
\]
If this is below $(n-6)/2$, no extension has defect at most five at that center. At a full clique, test its actual difference set against $D$. Every admissible set has a normalized representative surviving these reductions.

The program \texttt{exclude\_defect\_five.cpp} exhausts all four inputs without a time cutoff and finds no admissible set. The driver \texttt{verify\_defect\_five.py} regenerates the analytically complete inputs, checks their hashes and exact search counts, and verifies the analytic reductions and witness. Small direct subset enumerations test the pruning implementation, including a separate generating example with no unit difference. A separate full-modulus implementation checks the $\Z_{39}$ case. These checks are not an independent implementation of the entire large-modulus search, nor formal proof-assistant verification.
\end{proof}

\begin{table}[ht]
\centering
\begin{tabular}{@{}rrrrr@{}}\toprule
$n$ & Modulus bound & $|D|$ bound & Candidate $D$ & Recursive search nodes\\\midrule
16&685,700&41&550,496&630,696,891\\
18&4,114,202&39&349,351&282,562,384\\
20&24,685,212&39&349,167&216,055,373\\
22&39 only&39&1&15,708,790\\\bottomrule
\end{tabular}
\caption{Complete defect-five absence searches after the universal hole-budget filter. Graph and normalization reductions precede the recursive search. Every run terminates with no admissible configuration.}
\label{sup:tab:F5checks}
\end{table}

\begin{suptheorem}[Sharp threshold at defect five]\label{sup:thm:F5}
One has
\[
 \boxed{F(5)=14.}
\]
An extremizer is the complement in $\Z_{21}$ of
\[
 B=\{0,1,3,7,11,15,19\}.
\]
It has $(n,t,u,k,K)=(14,5,7,9,10)$. Every set of actual defect five with more than fourteen points has exactly five holes in its minimal regular completion.
\end{suptheorem}
\begin{proof}
Parity makes $n$ even. Theorem~\ref{res:2.1} excludes every $n>30$. For $n=24,26,28,30$, the hull bound gives $u\leq17$. Every odd $u\in\{7,9,11,13,15,17\}$ violates $(u-5)n\leq u(u-1)$ already at $n=24$. At $n=22$, the same hull bound gives $M\leq40$ and hence $u\leq17$. The budget excludes $u=7,9,11,13,15$, leaving $u=17$, $M=39$. Lemma~\ref{sup:lem:F5absence} excludes this case and the remaining sizes $16,18,20$. Therefore $F(5)\leq14$.

The displayed seven-point set is a classical sum cover \cite[Table~3]{haanpaa}; its covering property can also be checked directly. Its complement has fourteen points and full difference set. Since $4\cdot21>7^2+3\cdot7$, Proposition~\ref{res:5.2} gives actual defect $7-2=5$. Its hole count is seven, so exact completion fails. The pinned and global distance formulas give $k=9$, $K=10$. This proves the lower bound and the theorem.
\end{proof}

The construction itself is established in the sum-cover literature. The additional finite statement here is the exclusion of every larger counterexample over all odd cyclic completions. A classification of all fourteen-point extremizers is not asserted, and values $F(t)$ for $t\geq6$ remain unevaluated.

\section{Exceptional centers and effective recovery}\label{sup:sec:robust}
\begin{suptheorem}[Stability with exceptional centers]\label{sup:thm:robust}
Let $A$ be a nonempty subset of a finite abelian group of odd order, and let $s\geq0$ be an integer. Define
\[
 G=\{a\in A:\sig_A(a)\leq s\},\qquad
 b=|A\setminus G|,\qquad \tau=s+b.
\]
If $|G|>Q(\tau)$, let $C$ be the smallest affine subgroup coset containing $G$, and put $e=|A\setminus C|$. Then
\begin{equation}\label{sup:eq:robust}
 e\leq b,\qquad |C\setminus A|\leq s+e,\qquad
 |A\del C|\leq s+2e\leq s+2b.
\end{equation}
\end{suptheorem}
\begin{proof}
Removing $b$ points increases the reflection defect at a remaining center by at most $b$. More explicitly, for $a\in G$, the points $x\in G$ with $R_a(x)\notin G$ either satisfy $R_a(x)\notin A$, or are mapped injectively into $A\setminus G$. Thus
\[
 \sig_G(a)\leq\sig_A(a)+b\leq s+b=\tau.
\]
Let $t_G=\max_{a\in G}\sig_G(a)$. Since $Q$ is increasing on nonnegative integers and $|G|>Q(\tau)\geq Q(t_G)$, Theorem~\ref{res:2.1} gives $|C\setminus G|=t_G\leq\tau$. All points outside $C$ belong to $A\setminus G$, so $e\leq b$. Of the $b$ remaining points, exactly $b-e$ lie in $C\setminus G$. Therefore
\[
 |C\setminus A|=|C\setminus G|-(b-e)\leq(s+b)-(b-e)=s+e.
\]
Adding $e$ proves the symmetric-difference bounds.
\end{proof}

The theorem measures combinatorial deletions and outliers, not errors in Euclidean coordinates. No optimality claim is made for the coefficient of $b$.

For labeled vertices in $\Z_m$, the completion has a direct certificate. Compute the defects by testing whether $2a-x\pmod m$ belongs to $A$, for each ordered pair $(a,x)\in A^2$. Choose the good set $G$ when using Theorem~\ref{sup:thm:robust}; otherwise use $G=A$. With $a_0\in G$, put
\[
 g=\gcd\bigl(m,\{a-a_0:a\in G\}\bigr).
\]
Then $C=a_0+g\Z_m$ and $|C|=m/g$. This is the minimal completion, whether or not the theorem's size condition holds. The size condition certifies the stated bounds on its holes and outliers. A hash table gives $O(n^2)$ expected membership operations for the defects, followed by $O(n)$ gcd computations and $O(|C|)$ operations to enumerate the completion. These are operation counts; bit complexity also depends on $\log m$. The input consists of the vertex labels and $m$. Distance counts alone do not specify the location of the completion.

\section{Additional examples}\label{sup:sec:exceptions}
\begin{suptheorem}[An infinite local-to-global obstruction]\label{sup:thm:construction}
Let $q,r\geq3$ be odd integers, set $m=qr$, and define
\[
 B=\{0,1,\ldots,q-1\}\cup\{0,q,2q,\ldots,(r-1)q\},\qquad
 A=\Z_m\setminus B.
\]
Then the minimal affine completion is $\Z_m$ and
\begin{align*}
 n&=(q-1)(r-1),&u&=q+r-1,&t&=q+r-3,\\
 k&=(qr-3)/2,&K(A)&=(qr-1)/2=k+1.&
\end{align*}
In particular, for $q=r\to\infty$ through odd integers, $u=t+2$ and $n/t^2\to1/4$. There is no universal sufficient condition $n>f(t)$ for exact completion $u=t$ with $f(t)=o(t^2)$.
\end{suptheorem}
\begin{proof}
The two parts of $B$ intersect only at $0$, so $u=q+r-1$ and $n=qr-u=(q-1)(r-1)$. Each residue $s\in\Z_{qr}$ has a unique mixed-radix expression $s=j+iq$ with $0\leq j<q$ and $0\leq i<r$. Both summands belong to $B$. If $s\neq0$, these summands are distinct, since the two parts of $B$ intersect only at $0$. For $a\in A$, $2a\neq0$, so $r_{B+B}(2a)\geq2$. By \eqref{eq:5}, $\sig_A(a)\leq u-2$.

Take $a=q+1\in A$. We have $2a=2q+2<qr$. This residue cannot be a sum of two elements of $\{0,\ldots,q-1\}$, whose integer sums are at most $2q-2$. It cannot be a sum of two multiples of $q$ because $q\nmid2$. A sum with one element from each part is unique before reversing its order: it is $2+2q$. Such mixed sums lie between $0$ and $qr-1$, so no additional modular wraparound representation is possible. Thus $r_{B+B}(2a)=2$, and $t=u-2$.

The consecutive vertices $q+1,q+2$ belong to $A$, so its differences generate $\Z_m$. To prove that every difference occurs, observe that $2n>m$ except when $(q,r)=(3,3)$: the inequality is equivalent to $(q-2)(r-2)>2$. Density then gives $A-A=\Z_m$. In the remaining case, $A=\{4,5,7,8\}\subseteq\Z_9$, and its differences include $0,\pm1,\pm2,\pm3,\pm4$. Hence $A-A=\Z_m$ in all cases. Lemma~\ref{res:3.1} gives $k=(n-1+t)/2=(m-3)/2$, while the full difference set gives $K(A)=(m-1)/2$.

For $q=r$, $n=(q-1)^2$ and $t=2q-3$, so $n/t^2\to1/4$. If $f(t)=o(t^2)$, these examples eventually satisfy $n>f(t)$ but still have $u=t+2$, proving the final assertion.
\end{proof}

\subsection{Exact completion without symmetry}

\begin{supproposition}[Exact symmetry need not survive]\label{sup:prop:asymmetry}
The set $A=\Z_{31}\setminus\{0,1,3\}$ satisfies $n=28>Q(3)=18$ and $t=3$, but has no nonidentity rotation and no reflection preserving it.
\end{supproposition}
\begin{proof}
The defect assertion follows from Corollary~\ref{res:4.3}. A nonzero translation of $\Z_{31}$ has a single orbit of length $31$, so cannot preserve the three-element complement. If a reflection $x\mapsto c-x$ preserved that complement, its action on three points would have a fixed point $b\in\{0,1,3\}$; thus $c=2b\in\{0,2,6\}$. For $c=0$ the point $1$ maps to $30$; for $c=2$ the point $3$ maps to $30$; for $c=6$ the point $0$ maps to $6$. None preserves the complement, and therefore none preserves $A$.
\end{proof}

\begin{figure}[ht]
\centering\includegraphics[width=.94\linewidth]{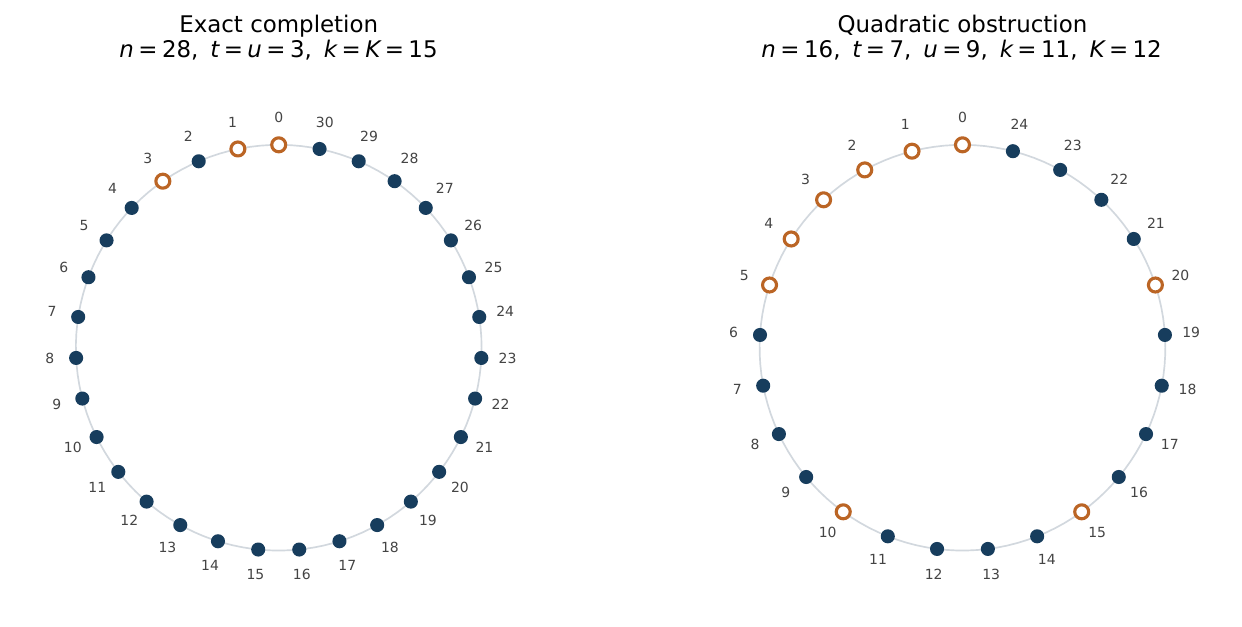}
\caption{Filled points belong to $A$; open points are holes in the ambient polygon. Left: exact completion with $t=3$, despite no dihedral symmetry of $A$. Right: the construction with $q=r=5$, for which $u=9>t=7$ and the global distance count exceeds the local maximum by one.}
\label{sup:fig:examples}
\end{figure}

\section{Exact computational checks}\label{sup:sec:checks}
The accompanying Python program uses integer modular arithmetic, with no floating-point distance comparisons. It computes pinned distances independently from unsigned cyclic steps, and reflection defects by membership tests. It checks \eqref{eq:3}, subgroup closure, the complement identity, the quadratic bound, and exact local-to-global completion whenever their hypotheses hold. Table~\ref{sup:tab:checks} reports the exhaustive runs. Normalization means that the set contains $0$; every nonempty set has a translate of this form, although this does not select one representative per translation orbit.

\begin{table}[ht]
\centering
\begin{tabular}{@{}rrr@{}}\toprule
Odd ambient order $m$ & Subsets tested & Exact-completion hypotheses met\\\midrule
3&4&1\\5&16&1\\7&64&1\\9&256&10\\11&1,024&11\\
13&4,096&13\\15&16,384&108\\17&65,536&137\\19&262,144&172\\\midrule
Total&349,524&454\\\bottomrule
\end{tabular}
\caption{Exhaustive nonempty subsets containing $0$. All checks passed. Singleton sets are included as boundary checks, although Corollary~\ref{res:2.2} is stated for $n\geq2$.}
\label{sup:tab:checks}
\end{table}

The program additionally checks all threshold choices for the exceptional-center statement in odd groups of order at most $13$, and $60$ structured examples with holes and outliers in larger cyclic groups. Eight members of Theorem~\ref{sup:thm:construction}, up to $(q,r)=(101,101)$, are checked using exact sum-representation counts. The example in Proposition~\ref{sup:prop:asymmetry} is checked against all rotations and reflections. The even-order identity in the main manuscript is tested on $10{,}922$ normalized subsets of groups of orders $2,4,\ldots,14$.

The retained program also tests the pointwise complement identity and the sum-cover equivalence for every one of the $349{,}524$ odd-order sets, and tests \eqref{eq:10} whenever all completion differences occur. It checks the representation budget, the finite hierarchy when its hypotheses hold, and the criterion \eqref{eq:12}. For the three-fibre family it checks $q\in\{7,13,19,31,61,121,301\}$, using an explicit Chinese-remainder labeling and exact ordered sum counts. The last example has $M=269{,}997$, $n=269{,}096$, $u=901$ and $t=899$. For $q=7,13$, pinned distances and differences are independently enumerated from $A$ itself.

The checks in this section are safeguards against errors in signs, parity, and boundary cases. They do not replace the analytic proofs, certify priority, or constitute formal verification in a proof assistant. The exhaustive computations in Lemmas~\ref{sup:lem:smallclass}, \ref{sup:lem:extendedclass} and~\ref{sup:lem:F5absence} are explicit parts of their computer-assisted proofs, with proved global bounds. Each proof states its cross-checks and their scope; an independent full implementation is not claimed for every search. The file \texttt{verification\_results.json} records the exact run, and \texttt{verify\_polygon.py} reproduces it using the Python standard library.

\section{An eventual finite sum-cover identity}\label{sup:sec:eventual_supp}
\subsection{A lower bound at every sufficiently large defect}
Define the finite odd cyclic sum-cover extremum
\begin{equation}\label{sup:eq:Sodd}
 S_{\rm odd}(u)=\max\{M\text{ odd}:\exists B\subseteq\Z_M,\ |B|=u,\ B+B=\Z_M\}.
\end{equation}
This is the odd-modulus restriction of the classical maximum-order sum-cover problem; compare \cite{haanpaa}. It is finite by $M\leq u(u+1)/2$.

\begin{supproposition}[Padding the three-fibre covers]\label{sup:prop:padding}
Let $u\geq158$, put $t=u-2$, and set
\[
 q=1+6\left\lfloor\frac{u-1}{18}\right\rfloor,\qquad M_0=3q(q-2),\qquad n_0=M_0-u.
\]
There is a cyclic set with $n_0$ points, actual defect $t$, $u$ holes, and $K-k=1$. Moreover,
\begin{align}
 n_0&>\max\left\{6t,\frac{(t+4)(t+3)}4\right\},\label{sup:eq:paddedregion}\\
 F(t)&\geq n_0\geq\frac{t^2-35t+241}{3}.\label{sup:eq:uniformlower}
\end{align}
\end{supproposition}
\begin{proof}
The three-fibre cover from Theorem~\ref{res:7.1} has $u_0=3q-2$ points in $\Z_{M_0}$. The definition of $q$ gives $0\leq u-u_0\leq17$. Add arbitrary distinct points to the cover until it has size $u$; the sum-cover property is preserved.

Since $q\geq(u-15)/3$ and $q\geq49$,
\[
 M_0\geq\frac{(u-15)(u-21)}3,\qquad
 n_0\geq\frac{u^2}{3}-13u+105.
\]
For $u\geq158$, this last quantity exceeds both $6(u-2)$ and $(u+2)(u+1)/4$. For the second comparison the difference is
\[
 \frac{u^2-165u+1254}{12}>0;
\]
it is positive at $158$ and increasing thereafter. The first comparison follows in the same way from $(u^2-57u+351)/3>0$. Also $M_0>2u$, so the complementary set has full difference set and minimal completion $\Z_{M_0}$. The second comparison implies $4M_0>u^2+3u$. Proposition~\ref{res:5.2} now forces $t=u-2$, and the gap identity gives $K-k=1$. Substituting $u=t+2$ gives \eqref{sup:eq:uniformlower}.
\end{proof}

\subsection{An exact reduction for every large fixed defect}
\begin{suptheorem}[Eventual finite sum-cover identity]\label{sup:thm:eventual}
For every integer $t\geq156$,
\begin{equation}\label{sup:eq:eventual}
 \boxed{F(t)=S_{\rm odd}(t+2)-(t+2).}
\end{equation}
Every extremal counterexample at such a defect has exactly $t+2$ holes and $K-k=1$. Its complement is a $(t+2)$-element sum cover of maximum possible odd cyclic order. Conversely, the complement of every such maximum-order cover is extremal for $F(t)$.
\end{suptheorem}
\begin{proof}
Put $u=t+2$. Proposition~\ref{sup:prop:padding} provides a counterexample of size $n_0$ strictly above the one-gap threshold in \eqref{sup:eq:paddedregion}. Thus any counterexample of maximum size $F(t)$ also lies above that threshold. Corollary~\ref{res:5.4} and failure of exactness force its number of holes to equal $t+2$. Its complement is a sum cover by Theorem~\ref{res:5.1}, so its completion order is at most $S_{\rm odd}(t+2)$. This proves the upper bound in \eqref{sup:eq:eventual}.

Conversely, let a cover of size $u$ attain $M=S_{\rm odd}(u)$. Then $M\geq M_0$, where $M_0$ is the padded construction's order. In particular $M>2u$ and $4M>u^2+3u$. Its complement has full difference set, actual defect $u-2=t$, and size $M-u$, by the same representation-count argument. It therefore supplies the reverse inequality. Both arguments also prove the extremal correspondence and the one-distance gap.
\end{proof}

The cutoff $156$ is explicit and sufficient; it is not asserted to be smallest. Equation~\eqref{sup:eq:eventual} gives a finite equality, stronger than an equality of limsups, but it still requires evaluation of the independent function $S_{\rm odd}$. It must not be read as a numerical formula for all $F(t)$.

\section{Amplification of strict sum covers}\label{sup:sec:amplification}
We now obtain the final lower bound while keeping every completion order odd. We use a finite strict sum cover as a seed and give an explicit amplification in cyclic groups. The numerical seed is classical; the transfer argument is proved here in full. It uses the two-direction covering mechanism of Bevan, Erskine and Lewis \cite[Lemma~5]{bevan}, specialized below to an elementary congruence calculation.

\begin{suplemma}[Two directions on a rectangular torus]\label{sup:lem:twodirections}
Let $C,d$ be positive integers with $d\mid C$, put $m=C/d$, and suppose $q>2C$, $r=q-C$, and $\gcd(q,C)=1$. Every $(x,y)\in\Z_q\times\Z_r$ can be written
\[
 (x,y)=h(1,1)+\ell(1,1+d),\qquad
 0\leq h<q+m,\quad 0\leq\ell<r.
\]
\end{suplemma}
\begin{proof}
Choose $0\leq x<q$, $0\leq y<r$. Since $d\mid C$ and $\gcd(r,C)=1$, let
\[
 \ell_0\equiv (y-x)d^{-1}\pmod r,\qquad 0\leq\ell_0<r.
\]
If $\ell_0\leq x$, take $\ell=\ell_0$ and $h=x-\ell_0$. Otherwise set
\[
 \ell=(\ell_0-m)\bmod r,\qquad h=x-\ell+q.
\]
In the latter case $h+\ell=x+q$ and
\[
 h+(1+d)\ell\equiv x+C+d(\ell_0-m)\equiv y\pmod r.
\]
If $\ell_0\geq m$, then $h=q+x-\ell_0+m<q+m$. If $\ell_0<m$, then $\ell=\ell_0-m+r$ and $h=C+x-\ell_0+m<C+m<q+m$. In both cases $h>0$, since $\ell<r<q$ and $x\geq0$. The first case has $0\leq h<q$ and is immediate. This proves all bounds and congruences.
\end{proof}

\begin{suptheorem}[Odd cyclic amplification]\label{sup:thm:amplification}
Let $w$ be odd and let $W=\{b_0,\ldots,b_{k-1}\}\subseteq\Z_w$ be a strict sum cover: every residue is $b_i+b_j$ for some $i<j$. Set
\[
 C=\operatorname{lcm}\bigl(2,w,2^1-1,2^2-1,\ldots,2^{k-1}-1\bigr).
\]
For every $q\equiv1\pmod C$ with $q>2C$, put $r=q-C$ and $L=q+C$. The set
\begin{equation}\label{sup:eq:amplifiedset}
 B_q=\bigcup_{i=0}^{k-1}
 \{(h\bmod q,\,2^ih\bmod r,\,b_i):0\leq h<L\}
 \subseteq\Z_q\times\Z_r\times\Z_w
\end{equation}
is a strict sum cover of a cyclic group of odd order
\[
 M_q=wq(q-C),\qquad |B_q|=k(q+C).
\]
Consequently,
\begin{equation}\label{sup:eq:seedlower}
 \liminf_{t\to\infty}\frac{F(t)}{t^2}\geq\frac{w}{k^2}.
\end{equation}
\end{suptheorem}
\begin{proof}
Since $C$ is even and divisible by $w$, the three orders $q,r,w$ are odd and pairwise coprime. Thus their product group is cyclic. Each direction $(1,2^i)$ has order $qr$, because $r$ is odd. Also $L<qr$ when $q>2C$ and $C\geq2$. Each fibre in \eqref{sup:eq:amplifiedset} therefore has $L$ distinct points, and its third coordinate separates it from the other fibres.

Given $(x,y,z)$, choose $i<j$ with $b_i+b_j=z$. Multiplying the second coordinate by the inverse of $2^i$ changes the two directions to $(1,1)$ and $(1,2^{j-i})$. The integer $d=2^{j-i}-1$ divides $C$. Lemma~\ref{sup:lem:twodirections} supplies the desired decomposition with parameters less than $q+C=L$ in the first fibre and less than $r<L$ in the second. The two summands lie in distinct fibres, so this is a strict representation.

For each sufficiently large integer $u$, choose the largest $q\equiv1\pmod C$ satisfying $k(q+C)\leq u$. Then
\[
 0\leq u-k(q+C)<kC,\qquad q=\frac{u}{k}+O_{W}(1).
\]
Add arbitrary distinct elements to $B_q$ to obtain exactly $u$ points; this is possible for large $u$ because $M_q$ grows quadratically. Covering is preserved. Hence
\[
 S_{\rm odd}(u)\geq M_q=\frac{w}{k^2}u^2-O_W(u).
\]
Apply Theorem~\ref{sup:thm:eventual} with $u=t+2$ to obtain \eqref{sup:eq:seedlower}. This is a uniform lower bound over all sufficiently large defects, not merely a subsequence bound.
\end{proof}

\begin{supcorollary}[An improved numerical lower bound]\label{sup:cor:improvedlower}
The optimal asymptotic completion coefficient satisfies
\begin{equation}\label{sup:eq:improvedlower}
 \boxed{\frac{53}{144}\leq
 \liminf_{t\to\infty}\frac{F(t)}{t^2}
 \leq\alpha_{\rm poly}=\beta_{\rm odd}\leq\frac12.}
\end{equation}
In particular $53/144=0.368055\ldots>1/3$.
\end{supcorollary}
\begin{proof}
Use the classical twelve-element strict sum cover
\[
 W=\{0,1,2,3,6,11,18,25,31,38,45,50\}\subseteq\Z_{53},
\]
listed in Fitch's tables accompanying \cite{fitch,fitchtables}. Directly listing its $\binom{12}{2}=66$ strict sums gives all $53$ residues. The accompanying certificate records a pair for each residue. Theorem~\ref{sup:thm:amplification} gives the lower bound; the other relations follow from Theorem~\ref{res:7.2} and the definition of $F$.
\end{proof}

The seed and its minimal-cover status are not discoveries of this paper; only its directly checked covering property is needed. The amplification argument is closely related to \cite{bevan}, and no claim of a best known constant across the literature is made. The explicit lower bound retained in this manuscript is $53/144$.

The construction has a substantial quantitative cost: with the displayed choice of slopes its fixed spacing is $C=3{,}673{,}365{,}502{,}978{,}710$. It is useful for the asymptotic question, not as an efficient way to evaluate small $F(t)$. The program \texttt{verify\_amplification.py} verifies the seed, exhaustive small two-direction cases, and exact decompositions at large moduli; the analytic proof covers every modulus in the stated family. Neither the equality $\alpha_{\rm poly}=1/2$ nor the existence of a limit $F(t)/t^2$ follows from this lower bound. All large-defect numerical evaluations still depend on determining $S_{\rm odd}$.

\section{Extension to odd-order finite abelian groups}
Define $\alpha_{\rm ab}$ by the formula for $\alpha_{\rm poly}$ with $A$ ranging over subsets of finite abelian groups of odd order, with $t$ the maximum reflection defect and $u$ measured in the minimal affine completion. Define $\beta_{\rm ab}$ by allowing all such groups $H$ in place of $\Z_M$, with numerator $|H|$.
\begin{suptheorem}
One has $\alpha_{\rm ab}=\beta_{\rm ab}$ and $1/3\le\alpha_{\rm ab}\le1/2$.
\end{suptheorem}
\begin{proof}
The decreasing suprema exist, and cyclic constructions give both lower bounds $1/3$. The completion theorem bounds the first constant by $1/2$; unordered sums bound the second. Choose counterexamples with $n/t^2\to\alpha_{\rm ab}$ and $t\to\infty$. Eventually $n>6t$, so the hull bound gives $u=O(t)$. The representation budget gives $0<u-t\le u(u-1)/n=O(1)$. Hence $u/t\to1$ and $(n+u)/u^2\to\alpha_{\rm ab}$. Their complements cover their minimal completion groups, proving $\alpha_{\rm ab}\le\beta_{\rm ab}$.
Conversely, take covers with $u\to\infty$ and $|H|/u^2\to\beta_{\rm ab}\ge1/3$. Their complements eventually have density greater than $1/2$, and therefore full difference set and minimal completion $H$. Since $4|H|>u^2+3u$ eventually, the representation budget forces $t=u-2$. Consequently $(|H|-u)/(u-2)^2\to\beta_{\rm ab}$, proving the reverse inequality.
\end{proof}
This extends the algebraic identification; it does not assign Euclidean polygon distances to noncyclic groups. The argument is a direct adaptation of the cyclic proof, and no separate priority claim is made. Numerical evaluation of these constants remains outside this result.

\FloatBarrier


\clearpage
\begin{thebibliography}{99}
\bibitem{bevan}
D. Bevan, G. Erskine and R. Lewis, Large circulant graphs of fixed diameter and arbitrary degree, \emph{Ars Math. Contemp.} \textbf{13} (2017), no. 2, 275--291. \href{https://arxiv.org/abs/1506.04962}{arXiv:1506.04962}.

\bibitem{chase}
Z. Chase, \emph{The maximum number of three term arithmetic progressions, and triangles in Cayley graphs}, preprint (2018), \href{https://arxiv.org/abs/1809.03729}{arXiv:1809.03729}.

\bibitem{custic}
A. \'{C}usti\'{c}, V. Kr\v{c}adinac and Y. Zhou, Tiling groups with difference sets, \emph{Electron. J. Combin.} \textbf{22} (2015), no. 2, Paper P2.56, 13 pp. \url{https://www.combinatorics.org/ojs/index.php/eljc/article/download/v22i2p56/pdf/0}.

\bibitem{delsarte}
P. Delsarte, J. M. Goethals and J. J. Seidel, Spherical codes and designs, \emph{Geom. Dedicata} \textbf{6} (1977), 363--388. \url{https://doi.org/10.1007/BF03187604}.

\bibitem{dumitrescu}
A. Dumitrescu, On distinct distances from a vertex of a convex polygon, \emph{Discrete Comput. Geom.} \textbf{36} (2006), 503--509. \url{https://doi.org/10.1007/s00454-006-1262-y}.

\bibitem{fitchtables}
M. A. Fitch, \emph{Sum covers: tables and example covers}, author-maintained companion tables; consulted 27 October 2026. \url{https://www.math.uaa.alaska.edu/~afmaf/math/sumcover/summary.html}.

\bibitem{fitch}
M. A. Fitch and R. E. Jamison, Minimum sum covers of small cyclic groups, \emph{Congr. Numer.} \textbf{147} (2000), 65--81.

\bibitem{green}
B. Green, \emph{100 open problems}, author-maintained problem collection, update December 2025, Problem~33, p.~18; consulted 27 October 2026. \url{https://people.maths.ox.ac.uk/greenbj/papers/open-problems.pdf}.

\bibitem{green-sisask}
B. Green and O. Sisask, On the maximal number of 3-term arithmetic progressions in subsets of \(\mathbb{Z}/p\mathbb{Z}\), \emph{Bull. Lond. Math. Soc.} \textbf{40} (2008), 945--955. \url{https://doi.org/10.1112/blms/bdn074}.

\bibitem{haanpaa}
H. Haanpää, Minimum sum and difference covers of Abelian groups, \emph{J. Integer Seq.} \textbf{7} (2004), Article 04.2.6. \url{https://cs.uwaterloo.ca/journals/JIS/VOL7/Haanpaa/haanpaa.pdf}.

\bibitem{jia}
X. Jia and J. Shen, Extremal bases for finite cyclic groups, \emph{SIAM J. Discrete Math.} \textbf{31} (2017), no. 2, 796--804. \url{https://doi.org/10.1137/15M1013237}.

\bibitem{kneser}
M. Kneser, Abschätzung der asymptotischen Dichte von Summenmengen, \emph{Math. Z.} \textbf{58} (1953), 459--484. \url{https://doi.org/10.1007/BF01174162}.

\bibitem{momihara}
K. Momihara and M. Shinohara, Distance sets on circles, \emph{Amer. Math. Monthly} \textbf{124} (2017), 241--254. \url{https://doi.org/10.4169/amer.math.monthly.124.3.241}.

\bibitem{nivasch}
G. Nivasch, J. Pach, R. Pinchasi and S. Zerbib, The number of distinct distances from a vertex of a convex polygon, \emph{J. Comput. Geom.} \textbf{4} (2013), 1--12. \href{https://arxiv.org/abs/1207.1266}{arXiv:1207.1266}.

\bibitem{nozaki}
H. Nozaki and M. Shinohara, On a generalization of distance sets, \emph{J. Combin. Theory Ser. A} \textbf{117} (2010), 810--826. \url{https://doi.org/10.1016/j.jcta.2009.11.001}.

\bibitem{shinohara}
M. Shinohara and K. Momihara, Distance sets and Kneser's theorem, \emph{RIMS Kôkyûroku} \textbf{2003} (2016), 96--103 (in Japanese). \url{https://www.kurims.kyoto-u.ac.jp/~kyodo/kokyuroku/contents/pdf/2003-12.pdf}.

\bibitem{tao}
T. Tao and V. Vu, \emph{Additive Combinatorics}, Cambridge Studies in Advanced Mathematics \textbf{105}, Cambridge University Press, 2006.
\end{thebibliography}
\end{document}